\documentclass[12pt,fleqn]{article}
\usepackage{latexsym,amsfonts,amssymb,amsmath,amsthm,longtable,epsfig,graphics,lscape}
\usepackage[cp1251]{inputenc}
\usepackage[english]{babel}

\newcommand{\ba} {\begin{array}}
\newcommand{\ea} {\end{array}}
\newcommand{\be} {\begin{equation}}
\newcommand{\ee} {\end{equation}}

\newtheorem{theorem}{Theorem}

{\theoremstyle{definition}

}

\begin{document}
\begin{center}
{\bf Dichotomy and bounded solutions of dynamical systems
in the Hilbert space.}
\end{center}

\begin{center}
Pokutnyi O.O. Institute of mathematics of NAS of Ukraine, lenasas@gmail.com
\end{center}
{\small Abstract. For a general discrete dynamics on a Banach and Hilbert spaces we give a necessary and
sufficient conditions of the existence of bounded solutions under assumption that the homogeneous difference equation
admits an exponential dichotomy on the semi-axes. We consider the so called resonance (critical) case when the uniqueness of solution is disturbed. We show that admissibility  can be reformulated in the terms of generalised or pseudoinvertibility. As an application we consider
the case when the corresponding dynamical system is e-trichotomy.    }

{\bf Introduction.}

{\bf Statement of the problem.}

Consider the following equation

\begin{equation} \label{6.2.1} x_{n+1} = A_{n}x_{n} + h_{n}, n \in
\mathbb{Z}, \end{equation}  where  $A_{n} : \mathbf{H} \rightarrow \mathbf{H} $ is a set of bounded operators which has a bounded inverse and
acts from the Hilbert space $H$ onto itself. Suppose that
$$ A = (A_{n})_{n \in
\mathbb{Z}} \in l_{\infty}(\mathbb{Z}, \mathcal{L}(H)), h =
(h_{n})_{n \in \mathbb{Z}} \in l_{\infty}(\mathbb{Z}, H).$$
It means that
$$|||A||| = \sup_{n \in \mathbb{Z}} ||A_{n}|| < + \infty, |||h||| =
\sup_{n \in \mathbb{Z}}||h_{n}|| < + \infty. $$ We need to establish conditions of the existence
of bounded solutions of the equation (\ref{6.2.1}).

Corresponding homogeneous equation has the following form:
\begin{equation} \label{6.2.2} x_{n + 1} = A_{n}x_{n}. \end{equation} It should be noted that any solution
of the homogeneous equation can be represented as: $x_{m} =
\Phi(m, n)x_{n}, m \geq n, $ where $\Phi(m, n) = A_{m
-1}A_{m-2}...A_{n+1}$ if $m >n$, and $\Phi(m, m) = I$. Clearly that $\Phi(m, 0) = A_{m -1}A_{m -
2}...A_{0},  U(m): = \Phi(m, 0), U(0) = I$.

Traditionally, the mapping $\Phi(m,n)$ is called evolution operator of the problem (\ref{6.2.2}).
Suppose that the equation (\ref{6.2.2}) is e-dichotomous on the semi axes  $\mathbb{Z}_{+}$ and $\mathbb{Z}_{-}$
with projectors $P$ and $Q$ in the space $H$ respectively, i.e.
$$
\exists k_{1} \geq 1; 0 < \lambda_{1} < 1; \exists P (P^{2} = P) :
$$
$$
||U(n)PU^{-1}(m)|| \leq k_{1}\lambda_{1}^{n - m}, n \geq m
$$
$$
||U(n)(E - P)U^{-1}(m)|| \leq k_{1}\lambda_{1}^{m - n}, m \geq n
$$
for any $m, n \in \mathbb{Z}_{+}$ (dichotomy on
$\mathbb{Z}_{+}$);
$$
\exists k_{2} \geq 1; 0 < \lambda_{2} < 1; \exists Q (Q^{2} = Q) :
$$
$$
||U(n)QU^{-1}(m)|| \leq k_{2}\lambda_{2}^{n - m}, n \geq m
$$
$$
||U(n)(E - Q)U^{-1}(m)|| \leq k_{2}\lambda_{2}^{m - n}, m \geq n
$$
for any $m, n \in \mathbb{Z}_{-}$ (dichotomy on
$\mathbb{Z}_{-}$).
\begin{theorem} \label{th6:2}  {\it Suppose that the homogeneous equation  is e-dichotomous on the semi-axes $\mathbb{Z}_{+}$ and $\mathbb{Z}_{-}$ with projectors $P$ and $Q$ respectively and the operator $D = P - (E - Q)$ is generalised invertible.  Solutions of the equation (\ref{6.2.1}) bounded on the whole axis exist if and only if the following condition is true
\begin{equation}
\label{6.2.3} \sum_{k = -\infty}^{+\infty}H(k+1)h_{k} = 0.
\end{equation}
Under condition (\ref{6.2.3}) the set of bounded solutions has the following form :
\begin{equation} \label{6.2.4}
x_{n}(c) = U(n)PP_{N(D)}c + (G[h])(n),
\end{equation}
where
$$
G[h](n) = U(n)Z(n),
$$
$$
Z(n) = \left\{
\begin{array}{rcl}
\sum_{k = 0}^{n - 1}PU^{-1}(k + 1)h_{k} - \sum_{k = n}^{+
\infty}(E - P)U^{-1}(k + 1)h_{k} + \\
+ PD^{-}[\sum_{k = 0}^{+\infty}(E - P)U^{-1}(k + 1)h_{k}+ \\
+ \sum_{k = -\infty}^{-1}QU^{-1}(k + 1)h_{k}],
\hspace{0.2cm} n \geq 0\\
\sum_{k = -\infty}^{n - 1}QU^{-1}(k + 1)h_{k} - \sum_{k =
n}^{-1}(E - Q)U^{-1}(k + 1)h_{k} + \\
+ (E - Q)D^{-}[\sum_{k = 0}^{+ \infty}(E - P)U^{-1}(k + 1)h_{k} +\\
+ \sum_{k = -\infty}^{-1}QU^{-1}(k + 1)h_{k}], \hspace{0.2cm} n
\leq 0,
\end{array}
\right.
$$
is generalised Green's operator on  $\mathbb{Z}$ with following properties :
$$
(G[h])(0 + 0) - (G[h])(0 - 0) = \sum_{k = -\infty}^{+ \infty} H(k + 1)h_{k} = 0,
$$
$$
(LG[h])(n) = h_{n}, n \in \mathbb{Z},
$$
where
$$
(Lx)(n):= x_{n + 1} - A_{n}x_{n} : l_{\infty}(\mathbb{Z}, B)
\rightarrow l_{\infty}(\mathbb{Z}, B),
$$
$H(n + 1) = P_{N(D^{*})}QU^{-1}(n + 1) = P_{N(D^{*})}(E -
P)U^{-1}(n + 1), D^{-}$  - is generalised invertible to operator $D$, $P_{N(D)}$ and  $P_{N(D^{*})}$  are projectors which project $\mathbf{B}$
onto kernel $N(D)$ and cokernel  $N(D^{*})$ of operators $D$ and
$D^{*}$ respectively.
}
\end{theorem}

{\bf Proof.} A general solution of the problem (\ref{6.2.1}), bounded on the semi axes has the form :
 \begin{equation} \label{6.2.5}
x_{n}(\xi) = \left\{
\begin{array}{rcl}
U(n)P\xi + \sum_{k = 0}^{n - 1}U(n)PU^{-1}(k + 1)h_{k} -\\
- \sum_{k = n}^{+ \infty}U(n)(E - P)U^{-1}(k + 1)h_{k},
\hspace{0.2cm} n \geq 0\\
U(n)(E - Q)\xi + \sum_{k = -\infty}^{n - 1}U(n)QU^{-1}(k + 1)h_{k}
- \\
- \sum_{k = n}^{-1}(E - Q)U^{-1}(k + 1)h_{k}, \hspace{0.2cm} n
\leq 0.
\end{array}
\right.
\end{equation} Prove that solution (\ref{6.2.5}) is bounded on the semi axes ($\mathbb{Z}_{+}, \mathbb{Z}_{-}$). Really, for any $n
\geq 0 $  we have  $A_{n}U(n)P\xi = A_{n}A_{n - 1}...A_{0}P\xi = U(n +
1)P\xi$. Then  $x_{n + 1} = A_{n}x_{n}$ and $x_{n}$ is solution of the homogeneous equation (\ref{6.2.2}) on $\mathbb{Z}_{+}$.
Further
$$
A_{n}(\sum_{k = 0}^{n - 1}U(n)PU^{-1}(k + 1)h_{k} - \sum_{k =
n}^{+\infty}U(n)(E - P)U^{-1}(k + 1)h_{k}) + h_{n} =
$$
$$
= \sum_{k = 0}^{n - 1}U(n + 1)PU^{-1}(k + 1)h_{k} - \sum_{k =
n}^{+\infty}U(n + 1)(E - P)U^{-1}(k + 1)h_{k} + h_{n} =
$$
$$
= \sum_{k = 0}^{n}U(n + 1)PU^{-1}(k + 1)h_{k} - \sum_{k = n +
1}^{+\infty}U(n + 1)(E - P)U^{-1}(k + 1)h_{k} + h_{n} -
$$
$$
 - U(n + 1)PU^{-1}(n + 1)h_{n} - U(n + 1)(E - P)U^{-1}(n +
1)h_{n} =
$$
$$
= x_{n + 1}(\xi). $$ Prove that the represented solution is bounded on the semi axes. Estimate series:
$$
||\sum_{k = 0}^{n - 1}U(n)PU^{-1}(k + 1)h_{k}|| \leq |||h|||
\sum_{k = 0}^{n - 1}||U(n)PU^{-1}(k + 1)|| \leq
$$
$$
\leq |||h|||\sum_{k = 0}^{n - 1}k_{1}\lambda_{1}^{n - k - 1} =
|||h|||k_{1}\lambda_{1}^{n}\sum_{k = 0}^{n - 1}\lambda_{1}^{-(k +
1)}=
$$
$$
=
|||h|||k_{1}\lambda_{1}^{n}\frac{\frac{1}{\lambda_{1}}((\frac{1}{\lambda_{1}})^{n
- 1} - 1)}{\frac{1}{\lambda_{1}} - 1} < \infty,
$$
and
$$ ||\sum_{k = n }^{+ \infty}U(n)(E - P)U^{-1}(k + 1)h_{k}||
\leq |||h|||\sum_{k = n}^{+\infty}k_{1}\lambda_{1}^{k + 1 - n} =
$$
$$
= k_{1}\lambda_{1}^{- n + 1} |||h|||\sum_{k =
n}^{+\infty}\lambda_{1}^{k} = k_{1}\lambda_{1}^{- n + 1}\cdot
\frac{\lambda_{1}^{n}}{1 - \lambda_{1}} < \infty .
$$
Boundedness of solution on  $\mathbb{Z}_{-}$ can be proved in such a way.

Find condition which guarantees that the solution  (\ref{6.2.5}) will be bounded on the whole integer axis. It will be if and only if when
$$
x_{0+}(\xi) = x_{0-}(\xi).
$$
Substitute corresponding expressions we obtain
$$
P\xi - \sum_{k = 0}^{+\infty}(E - P)U^{-1}(k + 1)h_{k} = (E -
Q)\xi + \sum_{k = - \infty }^{- 1}QU^{-1}(k + 1)h_{k}.
$$
Consider the following element
$$
g = \sum_{k = 0}^{+\infty}(E - P)U^{-1}(k + 1)h_{k} + \sum_{k = -
\infty}^{-1}QU^{-1}(k + 1)h_{k}.
$$
Obtain the following operator equation \begin{equation} \label{6.2.6} D\xi =
g.
\end{equation} Since $D$ is normally resolvable then as it is known \cite{BoiSam} a necessary and sufficient condition of solvability of the equation (\ref{6.2.6}) is following :
\begin{equation} \label{6.2.7} P_{N(D^{*})}g = 0.\end{equation}
Since $DP_{N(D)} = 0$, then we have that $PP_{N(D)} = (E -
Q)P_{N(D)}$. Since $P_{N(D^{*})}D = 0$ we have $P_{N(D^{*})}Q = P_{N(D^{*})}(E - P)$. Due to the equality condition
(\ref{6.2.7}) can be rewritten as
$$
\sum_{k = -\infty}^{ + \infty}P_{N(D^{*})}QU^{-1}(k + 1)h_{k} = 0,
$$
or
$$
\sum_{k = -\infty}^{+ \infty}P_{N(D^{*})}(E - P)U^{-1}(k + 1)h_{k}
= 0.
$$
In such a way we prove the condition (\ref{6.2.3}).
If the condition (\ref{6.2.3}) is true then $\xi = D^{-}g +
PP_{N(D)}c$, for any $c \in \mathbf{B}$.

Direct substitution in the representation (\ref{6.2.5}) gives us that the set of solutions bounded on the whole axis $\mathbb{Z}$ has the form (\ref{6.2.4}). $d$-normal and $n$-normal operators play important role in the theory of boundary value problems.

For such class of operators we can  obtain theorems with some refinements.

\begin{theorem}\label{th6:3} {\it Suppose that conditions of the theorem \ref{th6:2} are true and bounded operator $D = P - (E - Q)$
is $d$ - normal. Bounded on the whole axis solutions of the equation
 (\ref{6.2.1}) exist if and only if
$d$ linearly independent conditions are true \begin{equation} \label{6.2.8} \sum_{k
= -\infty}^{ + \infty} H_{d}(k + 1)h_{k} = 0. \end{equation} If the condition (\ref{6.2.8}) is true then solutions bounded on the whole axis have the  form (\ref{6.2.4}), where
$$
H_{d}(n) = [P_{N(D^{*})}Q]_{d}U^{-1}(n) = [P_{N(D^{*})}(E-
P)]_{d}U^{-1}(n),
$$
$$
d \leq m ~(m = dimcoker D < \infty), d = dim(P_{N(D^{*})}Q).
$$
}\end{theorem}

{\bf Proof.} It should be noted that the operator $P_{N(D^{*})}$ is finite dimensional (since  it is $d$ - normal), and then the operator $P_{N(D^{*})}Q$ is finite dimensional $(R(P_{N(D^{*})}Q)
\subset R(P_{N(D^{*})}))$.

\begin{theorem} \label{th6:4} {\it Suppose that conditions of the theorem \ref{th6:2} are true and the operator $D = P - (E - Q)$
is  $ n $ - normal. Then bounde on the whole axis solutions of the equation (\ref{6.2.1}) exist if and only if the condition (\ref{6.2.3}) is true.
Under condition (\ref{6.2.3}) the equation (\ref{6.2.1}) has  $r$ - parametric set of bounded solutions
\begin{equation} \label{6.2.9} x_{n}(c_{r}) =
U(n)[PP_{N(D)}]_{r}c_{r}  + (G[h])(n), \end{equation} where $r \leq
n$ ($n = dimker(D)$). }\end{theorem}

{\bf Proof.} Since  $D$ is $n$ - normal operator, then its kernel has finite dimension. It follows that the operator
$P_{N(D)}$ is finite dimensional and then the operator  $PP_{N(D)}$ is finite dimensional ($(R(PP_{N(D)}))\subset R(P_{N(D)})$).
If $dim~N(D) = n $, then
$$
dim(PP_{N(D)}B) = dim((E - Q)PP_{N(D)}) = r \leq n.
$$

\begin{theorem} \label{th6:5} {\it Supose that condition of the theorem \ref{th6:2} is hold and the operator $D = P -
(E - Q)$ is Fredholm. Then bounded solutions of the equation (\ref{6.2.1}) exist if and only if  $d$-conditions (\ref{6.2.8}) are true. Under condition (\ref{6.2.8}) the equation (\ref{6.2.1}) has
$r$-parametric set of bounded solutions
\begin{equation} \label{6.2.10} x_{n}(c_{r}) = U(n)[PP_{N(D)}]_{r}c_{r}
+ (G[h])(n), \end{equation} where $r \leq n$ ($n = dimker(D)$), $d
\leq m $ ($m = dimcoker D$). }\end{theorem}

{\bf Corollary.} {\it Suppose that under condition of the theorem \ref{th6:2} $[P,Q] =
PQ - QP = 0$ and $PQ = Q$. It is so called e-trichotomy case of the equation
(\ref{6.2.2}) on  $\mathbb{Z}$. In this case nonhomogeneous equation (\ref{6.2.1}) has at least one solution on
$\mathbb{Z}$ for any $h \in l_{\infty}(\mathbb{Z}, \mathbf{B})$. }

{\bf Proof.} From the equality $P_{N(D^{*})}D = 0$ and $DP = (P - (E -
Q))P = QP = Q$ follow that $P_{N(D^{*})}Q = P_{N(D^{*})}DP = 0$. From here we have solvability  $\forall h \in l_{\infty}(\mathbb{Z}, \mathbf{B})$.

{\bf Corollary 1.} {\it Under conditions of theorem (\ref{th6:2}) and additive condition  $[P,Q] =
PQ - QP = 0$ and $PQ = Q = P$ nonhomogeneous equation (\ref{6.2.1}) has unique bounded on $\mathbb{Z}$ solution for any
$ h \in l_{\infty}(\mathbb{Z}, \mathbf{B})$. }

{\bf Remark. } {\it In this case considering system is e-dichotomous on the whole axis $\mathbb{Z}$.
In the finite dimensional case the same result is well known \cite{Chuesh}. The theorem \ref{th6:2} under less restrictive assumptions
gives possibility to find the set of solutions. }

{\bf Corollary 2.} {\it Suppose that $[P, Q] = 0.$ Then operator $D$ has Moore-Penrose pseudoinvertible $D^{+}$ which is equal to $D$. In this case we have the following variants:
1 a.) equation (\ref{6.2.1}) has solutions if and only if the following condition is true ;

1 b.) under condition the set of bounded solutions has the form

2 a.) equation (\ref{6.2.1}) has bounded quasisolutions if and only if the following condition is true ;

2 b.) under condition the set of bounded solutions has the form

}
{\bf Proof.} If $[P, Q] = 0$ then
$$
DDD = (P - (I - Q))^{3} = P^{3} - 3 P^{2}(I - Q) + 3P(I - Q)^{2} - (I - Q)^{3} =
$$
$$
P - 3(P - PQ) + 3P(I - Q) - (I - Q) =  P - 3P + 3PQ + 3P - 3 PQ - I + Q = P - I + Q = D.
$$
From the equality we obtain that $D = D^{+}$.

In the Hilbert space we have more variants for the solutions.

\end{document}